\documentclass[12pt]{article}

\usepackage{amsmath}
\usepackage{amsthm}
\usepackage{amssymb}
\usepackage{times}
\usepackage[dvips, hyperindex]{hyperref}
\usepackage{enumerate}
\usepackage{amsfonts, latexsym}
\usepackage{stmaryrd}

\newtheorem{thm}{Theorem}
\newtheorem{lem}{Lemma}

\newtheorem{prop}{Proposition}

\newtheorem{question}{Question}

\newtheorem{cor}{Corollary}

\theoremstyle{definition}

\newcommand{\ZZ}{{\mathbf Z}}
\newcommand{\QQ}{{\mathbf Q}}

\newcommand{\CC}{{\mathbf C}}
\newcommand{\FF}{{\mathbb F}}

\newcommand{\OO}{{\mathcal O}}
\newcommand{\TT}{{\mathbb T}}
\newcommand{\cL}{{\mathcal L}}
\newcommand{\fP}{{\mathfrak P}}

\newcommand{\cross}{{\times}}

\newcommand{\Qbar}{{\overline{\QQ}}}
\newcommand{\Fbar}{{\overline{\FF}}}
\newcommand{\cO}{\mathcal{O}}

\renewcommand{\l}{{\ell}}
\DeclareMathOperator{\an}{an}

\newcommand{\GL}{\operatorname {GL}}
\newcommand{\SL}{\operatorname {SL}}
\newcommand{\PGL}{\operatorname {PGL}}
\newcommand{\PSL}{\operatorname {PSL}}
\newcommand{\Gal}{\operatorname {Gal}}
\newcommand{\End}{\operatorname {End}}

\newcommand{\frob}{\operatorname {Frob}}

\newcommand{\rhobar}{{\overline{\rho}}}

\newcommand{\mat}[4]{
 \left(  \begin{smallmatrix} #1 & #2 \\ #3 & #4 \end{smallmatrix} \right)}

\begin{document}
\title{On the generation of the coefficient field of a newform by a single Hecke eigenvalue}
\author{Koopa Tak-Lun Koo\footnote{Department of Mathematics,
University of Washington, Seattle, Box 354350 WA 98195, USA; e-mail:
{\tt koopakoo@gmail.com}}\,\, and
        William Stein\footnote{Department of Mathematics,
University of Washington, Seattle, Box 354350 WA 98195, USA; e-mail:
{\tt wstein@math.washington.edu}}\,\, and
        Gabor Wiese\footnote{Institut f\"ur Experimentelle Mathematik, Universit\"at Duisburg-Essen,
Ellernstra{\ss}e 29, 45326 Essen, Germany; e-mail: {\tt gabor.wiese@uni-due.de}}}

\maketitle

\begin{abstract}
  Let $f$ be a non-CM newform of weight $k \ge 2$.
  Let $L$ be a subfield of the coefficient field of~$f$.
  We completely settle the question of the density of the set of
  primes $p$ such that the $p$-th coefficient of~$f$ generates the
  field~$L$. This density is determined by the inner twists
  of~$f$. As a particular case, we obtain that in the absence of
  non-trivial inner twists, the density is~$1$ for $L$ equal to the whole
  coefficient field.
  We also present some new data on reducibility of Hecke polynomials, which
  suggest questions for further investigation.
  \medskip

\end{abstract}

\section{Statement of the results}\label{secone}

The principal result of this paper is the following theorem. Its corollaries below
completely resolve the question of the density of the set of primes~$p$ such that the $p$-th
coefficient of~$f$ generates a given field.

\begin{thm}\label{main density}
  Let $f$ be a newform (i.e., a new normalized cuspidal Hecke eigenform)
  of weight $k \ge 2$, level $N$ and Dirichlet character~$\chi$ which
  does not have complex multiplication (CM, see \cite[p.~48]{Rib80}).
  Let $E_f = \QQ(a_n(f) \,:\,(n,N)=1)$ be the field of coefficients
  of~$f$ and $F_f = \QQ\left(\frac{a_n(f)^2}{\chi(n)} \,:\,(n,N)=1\right)$.

  The set
  $$\left\{p \; \textnormal{prime}: \QQ\left(\frac{a_p(f)^2}{\chi(p)}\right) = F_f \right\}$$
  has density $1$.
\end{thm}

  A twist of~$f$ by a Dirichlet character~$\epsilon$ is said to be {\em inner}
  if there exists a (necessarily unique) field automorphism $\sigma_\epsilon: E_f \to E_f$
  such that 
\begin{equation}\label{eqin}
  a_p (f \otimes \epsilon) = a_p(f) \epsilon(p) = \sigma_\epsilon (a_p(f))
\end{equation}
  for almost all primes~$p$. For a discussion
  of inner twists we refer the reader to \cite[\S3]{Rib80} and
  \cite[\S3]{Rib85}. Here we give several statements that will be needed for
  the sequel.
  The $\sigma_\epsilon$ belonging to the inner twists
  of~$f$ form an abelian subgroup~$\Gamma$ of the automorphism group of~$E_f$.
  The field $F_f$ is the subfield of~$E_f$ fixed by~$\Gamma$.
  It is well-known that the coefficient field~$E_f$ is either a CM field
  or totally real. In the former case, the formula
\begin{equation}\label{eqcc}
   \overline{a_p(f)} = \chi(p)^{-1} a_p(f),
\end{equation}
  which is easily derived from the behaviour of the Hecke operators under the
  Petersson scalar product,
  shows that $f$ has a non-trivial inner twist by $\chi^{-1}$ with $\sigma_{\chi^{-1}}$
  being complex conjugation.
  If $N$ is square free, $k=2$ and the Dirichlet character $\chi$
  of $f$ is the trivial character,
  then there are no nontrivial inner twists of $f$.

\begin{lem}
  The field $F_f$ is totally real and $\QQ(a_p(f))$
  contains $\frac{a_p(f)^2}{\chi(p)}$.
\end{lem}

\begin{proof}
  Equation~\ref{eqcc} gives $\frac{a_p(f)^2}{\chi(p)} = a_p(f) \overline{a_p(f)}$,
  whence $F_f$ is totally real. Since every subfield of a CM field is preserved
  by complex conjugation, $\QQ(a_p(f))$ contains $\overline{a_p(f)}$, thus it also
  contains $\frac{a_p(f)^2}{\chi(p)}$.
\end{proof}

We immediately obtain the following two results.

\begin{cor}\label{cor density}
  Let $f$ and $E_f$ be as in Theorem~\ref{main density}.
  If $f$ does not have any nontrivial inner twists
  (e.g.\ if $k=2$, $N$ is square free and $\chi$ is trivial), then the set
  $$\left\{p \; \textnormal{prime}: \QQ(a_p(f)) = E_f \right\}$$
  has density $1$.
\end{cor}

\begin{cor}\label{cor Ff}
  Let $f$ and $F_f$ be as in Theorem~\ref{main density}.
  The set
  $$\left\{p \; \textnormal{prime}: F_f \subseteq \QQ(a_p(f)) \right\}$$
  has density $1$.
\end{cor}

  To any subgroup $H$ of $\Gamma$, we associate a number
  field $K_H$ as follows. Consider the inner twists as
  characters of the absolute Galois group $\Gal(\Qbar/\QQ)$ and let
  $\epsilon_1,\dots,\epsilon_r$ be the inner twists such that
  $H = \{\sigma_{\epsilon_1},\dots,\sigma_{\epsilon_r}\}$.
  Let $K_H$ be the minimal number field on which all $\epsilon_i$
  for $1\le i \le r$ are trivial, i.e.\ the field such that its absolute
  Galois group is the kernel of the map
  $$\Gal(\Qbar/\QQ) \xrightarrow{\epsilon_1,\dots,\epsilon_r} \CC^\times \times \dots \times \CC^\times.$$
  We use this field to express the density of the set of primes~$p$ such that
  $a_p(f)$ is {\em contained} in a given subfield of the coefficient field.

\begin{cor}
  Let $f$, $E_f$ and $F_f$ be as in Theorem~\ref{main density}. Let $L$
  be any subfield of~$E_f$. Let $M_L$ be the set
  $$ \left\{p \; \textnormal{prime}: a_p(f) \in L \right\}.$$

\begin{enumerate}[(a)]
\item If $L$ does not contain $F_f$, then $M_L$ has density~$0$.

\item If $L$ contains $F_f$, then $L = E_f^H$ for some subgroup $H \subseteq \Gamma$
  and $M_L$ has density $1/[K_H : \QQ]$.
\end{enumerate}
\end{cor}

\begin{proof}
Suppose first that $L$ does not contain~$F_f$. Then $a_p(f) \in L$ implies
that $F_f$ is not a subfield of~$\QQ(a_p(f))$. Thus by Corollary~\ref{cor Ff},
$M_L$ is a subset of a set of density~$0$ and is consequently itself of density~$0$.
We now assume that $L=E_f^H$. Then we have
\begin{align*}
 M_L 
 & = \left\{p \; \textnormal{prime}: \sigma(a_p(f)) = a_p(f) \, \forall \sigma \in H \right\}\\
 & = \left\{p \; \textnormal{prime}: a_p(f) \epsilon_i(p) = a_p(f) \, \forall i \in \{1,\dots,r\}\right\}.
\end{align*}
Since the set of~$p$ with $a_p(f)=0$ has density~$0$ (see for instance
\cite{Serre}, p.~174), the density of $M_L$ is equal to the density of 
$$\left\{p \; \textnormal{prime}: \epsilon_i(p) = 1 \, \forall i \in \{1,\dots,r\} \right\}
= \left\{p \; \textnormal{prime}: p \textnormal{ splits completely in } K_H \right\},$$
yielding the claimed formula.
\end{proof}

  A complete answer as to the density of the set of~$p$ such that $a_p(f)$ {\em generates}
  a given field $L \subseteq E_f$ is given by the following immediate result.

\begin{cor}
  Let $f$, $E_f$ and $F_f$ be as in Theorem~\ref{main density}. Let $L$
  be $E_f^H$ with $H$ some subgroup of~$\Gamma$. The density of the set
  $$ \left\{p \; \textnormal{prime}: \QQ(a_p(f)) = L \right\}.$$
  is equal to the density of the set
  $$\left\{p \; \textnormal{prime}: \epsilon_i(p) = 1 \, \forall i \in \{1,\dots,r\} \textnormal{ and }
     \epsilon_j(p) \neq 1 \, \forall j \in \{r+1,\dots,s\} \right\},$$
  where the $\epsilon_j$ for $j \in \{r+1,\dots,s\}$ are the inner twists of~$f$ that
  belong to elements of~$\Gamma - H$.
\end{cor}

  This corollary means that the above density is completely
  determined by the inner twists of~$f$.
  We illustrate this by giving two examples.
  In weight~$2$ there is a newform on $\Gamma_0(63)$ with coefficient field~$\QQ(\sqrt{3})$.
  It has an inner twist by the Legendre symbol $p \mapsto \left(\frac{p}{3}\right)$.
  Consequently, the field $F_f$ is $\QQ$ and the set of~$p$ such that
  $a_p(f) \in \QQ$ has density~$\frac{1}{2}$.

  For the next example we consider the newform of weight~$2$ on $\Gamma_0(512)$
  whose coefficient field has degree~$4$ over~$\QQ$.
  More precisely, the coefficient field $E_f$ is $\QQ(\sqrt{2},\sqrt{3})$ and $F_f = \QQ$.
  Hence, $\Gamma = \ZZ/2\ZZ \times \ZZ/2\ZZ = \{1, \sigma_1,\sigma_2,\sigma_3\}$.
  There are thus nontrivial inner twists $\epsilon_1$, $\epsilon_2$ and $\epsilon_3$,
  all of which are quadratic, as
  their values must be contained in the totally real field~$E_f$.
  As $\sigma_1 \sigma_2 = \sigma_3$,
  it follows that $\epsilon_1(p) \epsilon_2(p) = \epsilon_3(p)$.
  This equation already excludes the possibility that all $\epsilon_i(p) \neq 1$,
  whence there is not a single $p$ such that $a_p(f)$ generates~$E_f$.
  Furthermore, the set of~$p$ such that $a_p$ generates the quadratic field
  $E_f^{\langle \sigma_1 \rangle}$ is equal to the density of
  $\left\{p \; \textnormal{prime}: \epsilon_1(p)=1 \textnormal{ and } \epsilon_2(p) \neq 1 \right\},$
  which is~$\frac{1}{4}$. Similar arguments apply to the other two quadratic fields.
  The set of~$p$ such that $a_p \in \QQ$ also has density~$\frac{1}{4}$.

  \medskip

  In the literature there are related but weaker results concerning
  Corollary~\ref{cor density}, which are situated
  in the context of Maeda's conjecture, i.e., they concern the case of
  level~$1$ and assume that the space $S_k(1)$ of cusp forms of weight~$k$
  and level~$1$ consists of a single Galois orbit
  of newforms (see, e.g., \cite{JOno} and \cite{BM2003}).
  We now show how Corollary~\ref{cor density} extends the
  principal results of these two papers.

  Let $f$ be a newform of level~$N$, weight~$k \ge 2$
  and trivial Dirichlet character $\chi=1$ which neither has CM nor 
  nontrivial inner twists. This is for instance true when $N=1$.
  Let $\TT$ be the $\QQ$-algebra generated by all $T_n$ with $n\ge 1$
  inside $\End(S_k(N,1))$ and let $\fP$ be the kernel of the
  $\QQ$-algebra homomorphism $\TT \xrightarrow{T_n \mapsto a_n(f)} E_f$.
  As $\TT$ is reduced, the map $\TT_\fP \xrightarrow{T_n \mapsto a_n(f)} E_f$
  is a ring isomorphism with $\TT_\fP$ the localization of $\TT$ at~$\fP$.
  Non canonically $\TT_\fP$ is also isomorphic as a $\TT_\fP$-module
  (equivalently as an $E_f$-vector space) to its $\QQ$-linear dual, which can be
  identified with the localization at~$\fP$ of the $\QQ$-vector space $S_k(N,1; \QQ)$
  of cusp forms in $S_k(N,1)$ with $q$-expansion in $\QQ[[q]]$.
  Hence, $\QQ(a_p(f)) = E_f$ precisely means that the characteristic
  polynomial~$P_p \in \QQ[X]$ of~$T_p$ acting on the localization at~$\fP$ 
  of $S_k(N,1;\QQ)$ is irreducible.
  Corollary~\ref{cor density} hence shows that the set of
  primes~$p$ such that $P_p$ is irreducible has density~$1$. 

  This extends Theorem 1 of~\cite{JOno} and Theorem 1.1 of~\cite{BM2003}.
  Both theorems restrict to the case $N=1$ and assume that there is a
  unique Galois orbit of newforms, i.e., a unique~$\fP$, so that no 
  localization is needed. Theorem~1 of~\cite{JOno} says that
  $$\# \{p < X \, \textnormal{prime} \, :
     P_p \textnormal{ is irreducible in } \QQ[X]\} \gg \frac{X}{\log X}$$
  and Theorem 1.1 of~\cite{BM2003} states that there is $\delta > 0$
  such that
  $$\# \{p < X \, \textnormal{prime} \, :
     P_p \textnormal{ is reducible in } \QQ[X]\} \ll \frac{X}{(\log X)^{1+\delta}}.$$

\vspace{2ex}
\noindent{\bf Acknowledgements.}
  The authors would like to thank the MSRI, where part of this research was
  done, for its hospitality. The first author would like to thank his
  advisor Ralph Greenberg for suggesting the problem. The second author
  acknowledges partial support from the National Science
  Foundation grant No.\ 0555776, and also used \cite{sage} for some calculations
  related to this paper.
  All three authors thank Jordi Quer for useful discussions.

\section{Group theoretic input}

\begin{lem}\label{fh}
  Let $q$ be a prime power and $\epsilon$ a generator of the cyclic
  group $\FF_q^{\cross}.$

\begin{enumerate}[(a)]

\item The conjugacy classes $c$ in $\GL_2(\FF_q)$ have the
  following four kinds of representatives:
$$   S_a   = \begin{pmatrix}a & 0 \\ 0 & a \end{pmatrix}, \; \;
     T_{a} = \begin{pmatrix}a & 0 \\ 1 & a \end{pmatrix}, \; \;
  U_{a, b} = \begin{pmatrix}a & 0 \\ 0 & b \end{pmatrix}, \; \;
  V_{x, y} = \begin{pmatrix}x & \epsilon y \\ y & x \end{pmatrix}$$
  where $a \not= b,$ and $y \not= 0.$

\item The number of elements in each of these conjugacy classes are:
  $1, q^2 - 1, q^2 + q,$ and $q^2 - q$, respectively.

\end{enumerate}
\end{lem}

\begin{proof} See Fulton-Harris~\cite{FH}, page 68.
\end{proof}

We use the notation $[g]_G$ for the conjugacy class of~$g$ in~$G$.

\begin{prop}\label{gpthy}
Let $q$ be a prime power and $r$ a positive integer.
Let further $R \subseteq \widetilde{R} \subseteq \FF_{q^{r}}^\times$ be subgroups.
Put $\sqrt{\widetilde{R}} = \{ s \in \FF_{q^r}^\times \,:\, s^2 \in \widetilde{R}\}$.
Set
$$
H = \{ g \in \GL_2(\FF_{q}) \,:\, \det (g) \in R \}
$$
and let
$$
 G \subseteq \{ g \in \GL_2(\FF_{q^r}) \,:\, \det (g) \in \widetilde{R} \}
$$
be any subgroup such that $H$ is a normal subgroup of~$G$.
Then the following statements hold.
\begin{enumerate}[(a)]

\item The group $G/(G \cap \FF_{q^r}^\times)$ (with $\FF_{q^r}^\times$ identified
with scalar matrices) is either equal to $\PSL_2(\FF_q)$
or to $\PGL_2(\FF_q)$. More precisely, if we let $\{s_1,\dots,s_n\}$ be a system
of representatives for $\sqrt{\widetilde{R}}/R$, then for all $g \in G$
there is $i$ such that $g \mat {s_i^{-1}} 00 {s_i^{-1}} \in G \cap \GL_2(\FF_q)$
and $\mat {s_i}00{s_i} \in G$.

\item Let $g \in G$ such that $g \mat {s_i^{-1}} 00 {s_i^{-1}} \in G \cap \GL_2(\FF_q)$
and $\mat {s_i}00{s_i} \in G$. Then
$$[g]_G = [g \mat {s_i^{-1}} 00 {s_i^{-1}}]_{G \cap \GL_2(\FF_q)} \mat {s_i}00{s_i}.$$

\item Let $P(X) = X^2 - aX + b \in \FF_{q^r}[X]$ be a polynomial.
  Then the inequality
  $$ \sum_C | C | \;\le\; 2 | \widetilde{R}/R | (q^2 + q)$$
  holds, where the sum runs over the conjugacy classes $C$ of~$G$
  with characteristic polynomial equal to~$P(X)$.

\end{enumerate}

\end{prop}

\begin{proof}
(a) The classification of the finite subgroups of $\PGL_2(\Fbar_q)$
yields that the group
$G/(G \cap \FF_{q^r}^\times)$ is either $\PGL_2(\FF_{q^u})$ or $\PSL_2(\FF_{q^u})$
for some $u \mid r$.
This, however, can only occur with $u=1$, as
$\PSL_2(\FF_{q^u})$ is simple.
The rest is only a reformulation.

(b) This follows from (a), since scalar matrices are central.

(c) From (b) we get the inclusion
$$\bigsqcup_C C \subseteq \bigsqcup_{i=1}^n \bigsqcup_D D \mat {s_i}00{s_i},$$
where $C$ runs over the conjugacy classes of $G$ with characteristic
polynomial equal to $P(X)$ and $D$ runs over the conjugacy classes
of $G\cap \GL_2(\FF_q)$ with characteristic polynomial equal to $X^2
- a s_i^{-1} X + b s_i^{-2}$ (such a conjugacy class is empty if the
polynomial is not in $\FF_q[X]$). The group $G\cap \GL_2(\FF_q)$ is
normal in $\GL_2(\FF_q)$, as it contains $\SL_2(\FF_q)$. Hence, any
conjugacy class of $\GL_2(\FF_q)$ either has an empty intersection with
$G\cap \GL_2(\FF_q)$ or is a disjoint union of conjugacy classes of
$G\cap \GL_2(\FF_q)$. Consequently, by Lemma~\ref{fh}, the disjoint union
$\bigsqcup_D D \mat {s_i}00{s_i}$ is equal to one of
\begin{enumerate}[(i)]
\item $[U_{a,b}]_{\GL_2(\FF_q)} \mat {s_i}00{s_i}$,
\item $[V_{x,y}]_{\GL_2(\FF_q)} \mat {s_i}00{s_i}$ or
\item $[S_a]_{\GL_2(\FF_q)} \mat {s_i}00{s_i} \sqcup  [T_a]_{\GL_2(\FF_q)} \mat {s_i}00{s_i}$.
\end{enumerate}
Still by Lemma~\ref{fh}, the first set contains $q^2 + q$, the
second set $q^2-q$ and the third one $q^2$ elements. Hence, the set
$\bigsqcup_C C$ contains at most $2 | \widetilde{R}/R | (q^2+q)$
elements.
\end{proof}

\section{Proof}

  The proof of Theorem~\ref{main density} relies on the following
  important theorem by Ribet, which, roughly speaking, says that the
  image of the mod~$\l$ Galois representation attached to a fixed
  newform is as big as it can be for almost all primes~$\l$.

\begin{thm}[Ribet]\label{ribet}
  Let $f$ be a Hecke eigenform of weight $k \ge 2$, level~$N$ and Dirichlet
  character $\chi: (\ZZ/N\ZZ)^\times \to \CC^\times$. Suppose that $f$
  does not have CM.
  Let $E_f$ and $F_f$ be as in Theorem~\ref{main density} and
  denote by $\cO_{E_f}$ and $\cO_{F_f}$ the corresponding rings
  of integers.
  For almost all prime numbers~$\l$ the following statement holds:
\begin{quote}
  Let $\widetilde{\cL}$ be a prime ideal of $\cO_{E_f}$ dividing~$\l$.
  Put $\cL = \widetilde{\cL} \cap \cO_{F_f}$ and $\cO_{F_f}/\cL \cong \FF$.
  Consider the residual Galois representation
  $$\rhobar_{f,\widetilde{\cL}}: \Gal(\Qbar/\QQ) \to \GL_2(\cO_{E_f}/\widetilde{\cL})$$
  attached to~$f$. Then the image $\rhobar_{f,\widetilde{\cL}}(\Gal(\Qbar/K_\Gamma))$ is equal to
  $$ \{g \in \GL_2(\FF) \,:\, \det(g) \in \FF_\l^{\times (k-1)} \},$$
  where $K_\Gamma$ is the field defined in Section~\ref{secone}.
\end{quote}
\end{thm}

\begin{proof}
  It suffices to take Ribet~\cite[Thm.~3.1]{Rib85} mod~$\widetilde{\cL}$.
\end{proof}

\begin{thm}\label{global_density}
  Let $f$ be a non-CM newform of weight $k \ge 2$, level $N$ and Dirichlet
  character~$\chi$.
  Let $F_f$ be as in Theorem~\ref{main density}
  and let $L \subset F_f$ be any proper subfield. Then the set
  $$\left\{p \; \textnormal{prime}: \frac{a_p(f)^2}{\chi(p)} \in L\right\}$$
  has density zero.
\end{thm}

\begin{proof}
  Let $L \subsetneq F_f$ be a proper subfield and $\OO_L$ its integer ring.
  We define the set
  $$ S := \{ \cL \subset \OO_{F_f} \textnormal{ prime ideal}:
  [\OO_{F_f}/\cL : \OO_L / (L \cap \cL)] \ge 2 \}.$$
  Notice that this set is infinite. For, if it were finite, then all
  but finitely many primes would split completely in the extension
  $F_f/L$, which is not the case by Chebotarev's density theorem.

  Let $\cL \in S$ be any prime, $\l$ its residue characteristic and
  $\widetilde{\cL}$ a prime of $\OO_{E_f}$ lying over~$\cL$.
  Put $\FF_q = \OO_L / (L \cap \cL)$,
      $\FF_{q^r} = \OO_{F_f}/\cL$ and
      $\FF_{q^{rs}} = \OO_{E_f}/\widetilde{\cL}$. We have $r \ge 2$.
  Let $W$ be the subgroup of $\FF_{q^{rs}}^\times$ consisting of the values
  of~$\chi$ modulo~$\widetilde{\cL}$; its size $|W|$ is less than or equal
  to $| (\ZZ/N\ZZ)^\times|$.
  Let
    $R = \FF_\l^{\times (k-1)}$
  be the subgroup of $(k-1)$st powers of elements in the multiplicative
  group $\FF_\l^{\times}$ and let
    $\widetilde{R} = \langle R, W \rangle \subset \FF_{q^{rs}}^\times$.
  The size of $\widetilde{R}$ is less than or equal to $|R| \cdot |W|$.
  Let
    $H = \{g \in \GL_2(\FF_{q^r}) \,:\, \det(g) \in R \}$
  and
    $G = \Gal(\Qbar^{\ker{\rhobar_{f,\widetilde{\cL}}}}/\QQ)$.
  By Galois theory, $G$ can be identified with the image of the residual representation
  $\rhobar_{f,\widetilde{\cL}}$, and we shall make this identification from
  now on.
  By Theorem~\ref{ribet} we have the inclusion of groups
  $$ H \subseteq G \subseteq  \{g \in \GL_2(\FF_{q^{rs}}): \det(g) \in \widetilde{R} \}$$
  with $H$ being normal in~$G$.

  If $C$ is a conjugacy class of~$G$, by Chebotarev's
  density theorem the density of
  $$\{p \, \textnormal{prime} : \, [\rhobar_{f,\widetilde{\cL}}(\frob_p)]_G = C\}$$
  equals $|C|/|G|.$
  We consider the set
  $$ M_\cL := \bigsqcup_C \{ p \, \textnormal{prime} : \, [\rhobar_{f,\widetilde{\cL}}(\frob_p)]_G = C\}
     \supseteq 
   \left\{ p \, \textnormal{prime} : \, \overline{\left(\frac{a_p(f)^2}{\chi(p)}\right)} \in \FF_q \right\},$$
  where the reduction modulo~$\cL$ of an element $x \in \cO_{F_f}$ is denoted by $\overline{x}$
  and $C$ runs over the conjugacy classes of~$G$ with
  characteristic polynomials equal to some $X^2-aX+b \in \FF_{q^{rs}}[X]$
  such that 
  $$a^2 \in \{ t \in \FF_{q^{rs}} \, : \, \exists u \in \FF_q \; \exists w \in W : t = uw \}$$
  and automatically $b \in  \widetilde{R}$.
  The set $M_\cL$ has the density $\delta(M_\cL) = \sum_{C}^{}\frac{|C|}{|G|}$
  with $C$ as before.
  There are at most $2q |W|^2 \cdot |R|$ such polynomials.
  We are now precisely in the situation to apply Prop.~\ref{gpthy}, Part~(c),
  which yields the inequality
  $$     \delta (M_\cL)
     \le \frac{4 |W|^3 q (q^{2r} +q^{r})}{(q^{3r} - q^r)}
     =   O\left(\frac{1}{q^{r - 1}}\right) \le O\left(\frac{1}{q}\right),                $$
  where for the denominator we used $|G| \ge |H| = |R| \cdot |\SL_2(\FF_{q^r})|$.

  Since $q$ is unbounded for $\cL \in S$,
  the intersection $M := \bigcap_{\cL \in S} M_\cL$ is a set having a density
  and this density is~$0$.
  The inclusion
  $$ \left\{ p \, \textnormal{prime} : \, \frac{a_p(f)^2}{\chi(p)} \in L \right\}
       \subseteq M $$
  finishes the proof.
\end{proof}

\begin{proof}[Proof of Theorem~\ref{main density}]
  It suffices to apply Theorem~\ref{global_density} to
  each of the finitely many subextension of~$F_f$.
\end{proof}

\section{Reducibility of Hecke polynomials: questions}

Motivated by a conjecture of Maeda, there has been some speculation
that for every integer $k$ and prime number $p$, the characteristic
polynomial of $T_p$ acting on $S_k(1)$ is irreducible.  See, for
example, \cite{fj}, which verifies this for all $k<2000$ and $p<2000$.
The most general such speculation might be the following question: 
{\em if~$f$ is a non-CM newform of level
$N\geq 1$ and weight $k \ge 2$ such that some $a_p(f)$ generates the field $E_f =
\QQ(a_n(f) : n\geq 1)$, do all but finitely many prime-indexed Fourier
coefficients $a_p(f)$ generate~$E_f$?}
The answer in general is no. An example is given by the newform in level~$63$
and weight~$2$ that has an inner twist by $\left(\frac{\cdot}{3}\right)$.
Also for non-CM newforms of weight~$2$ without nontrivial inner twists such that
$[E_f:\QQ]=2$, we think that the answer is likely no.

Let $f\in S_k(\Gamma_0(N))$ be a newform of weight $k$ and level $N$.
The {\em degree} of $f$ is the degree of the field $E_f$, and we say
that~$f$ is a {\em reducible newform} if $a_p(f)$ does not generate~$E_f$
for infinitely many primes~$p$.

For each even weight $k\leq 12$ and degree $d=2,3,4$, we used
\cite{sage} to find newforms $f$ of weight $k$ and degree $d$.  For
each of these forms, we computed the {\em reducible primes} $p<1000$,
i.e., the primes such $a_p(f)$ does not generate~$E_f$.
The result of this computation is given in Table~\ref{counting}.
Table~\ref{newforms2} contains the number of reducible primes
$p<10000$ for the first $20$ newforms of degree $2$ and weight $2$.
This data inspires the following question.
\begin{question}
If $f\in S_2(\Gamma_0(N))$ is a newform of degree $2$,
is $f$ necessarily reducible?  That is, are there
infinitely many primes $p$ such that $a_p(f)\in \ZZ$?
\end{question}

Tables~\ref{newforms3}--\ref{newforms389} contain additional
data about the first few newforms of given degree and weight, 
which may suggest other similar questions.
In particular, Table~\ref{tab:million} contains data for
all primes up to $10^6$ for the first degree 2 form $f$
with $L(f,1)\neq 0$, and for the first degree 2 form
$g$ with $L(g,1) = 0$.  We find that there are 386 primes
$<10^6$ with $a_p(f) \in \ZZ$ and $309$ with $a_p(g)\in \ZZ$.

\begin{question}
If $f\in S_2(\Gamma_0(N))$ is a newform of degree $2$,
can the asymptotic behaviour of the function
$$ N(x) := \#\{ p \, \textnormal{prime} : \, p < x, a_p(f) \in \ZZ \}  $$
be described as a function of~$x$?
\end{question}

The authors intend to investigate these questions in a subsequent paper.

\begin{table}
\begin{center}
\caption{Counting Reducible Characteristic Polynomials\label{counting}}
\vspace{0.5ex}
\begin{tabular}{|l|l|l|l|}\hline
$k$ & $d$ & $N$ & reducible $p<1000$\\\hline
2 & 2 & 23 & 13, 19, 23, 29, 43, 109, 223, 229, 271, 463, 673, 677, 883, 991\\
2 & 3 & 41 & 17, 41\\
2 & 4 & 47 & 47 \\\hline
4 & 2 & 11 & 11 \\
4 & 3 & 17 & 17 \\
4 & 4 & 23 & 23 \\\hline
6 & 2 & 7 & 7 \\
6 & 3 & 11 & 11 \\
6 & 4 & 17 & 17 \\\hline
8 & 2 & 5 & 5 \\
8 & 3 & 17 & 17 \\
8 & 4 & 11 & 11 \\\hline
10 & 2 & 5 & 5 \\
10 & 3 & 7 & 7 \\
10 & 4 & 13 & 13 \\\hline
12 & 2 & 5 & 5 \\
12 & 3 & 7 & 7 \\
12 & 4 & 21 & 3, 7 \\
\hline
\end{tabular}
\end{center}
\end{table}

\begin{table}
\begin{center}
\caption{First 20 Newforms of Degree 2 and Weight 2\label{newforms2}}
\vspace{0.5ex}
\begin{tabular}{|l|l|l|c|}\hline
$k$ & $d$ & $N$ & \#\{reducible $p<10000$\}\!\\\hline
2 & 2 & 23 & 47 \\
2 & 2 & 29 & 42 \\
2 & 2 & 31 & 78 \\
2 & 2 & 35 & 48 \\
2 & 2 & 39 & 71 \\
2 & 2 & 43 & 43 \\
2 & 2 & 51 & 64 \\
2 & 2 & 55 & 95 \\
2 & 2 & 62 & 77 \\
2 & 2 & 63 & 622 (inner twist by $\left(\frac{\cdot}{3}\right)$)\\
\hline
\end{tabular}\hspace{1em}\begin{tabular}{|l|l|l|c|}\hline
$k$ & $d$ & $N$ & \#\{reducible $p<10000$\}\!\\\hline
2 & 2 & 65 & 43 \\
2 & 2 & 65 & 90 \\
2 & 2 & 67 & 51 \\
2 & 2 & 67 & 19 \\
2 & 2 & 68 & 53 \\
2 & 2 & 69 & 47 \\
2 & 2 & 73 & 43 \\
2 & 2 & 73 & 55 \\
2 & 2 & 74 & 52 \\
2 & 2 & 74 & 21 \\
\hline
\end{tabular}
\end{center}
\end{table}

\begin{table}\label{tab:million}
\begin{center}
\caption{Newforms 23a and 67b: values of $\psi(x) = \#\{\text{reducible }p< x\cdot 10^5\}$
\label{newforms1000000}}
\vspace{0.5ex}
\begin{tabular}{|c|c|c|c|c|c|c|c|c|c|c|c|c|c|}\hline
$k$ & $d$ & $N$ & $r_{\an}$ & $1$ & $2$ & $3$ &$4$ &$5$ &$6$ &$7$ &$8$ &$9$ & $10$ \\\hline
2 & 2 & $23$ &  0 & 127 & 180 & 210 & 243 & 277 & 308 & 331 & 345 & 360 & 386
\\\hline
2 & 2 & $67$ & 1 & 111 & 159 & 195 & 218 & 240 & 257 & 276 & 288 & 301 & 309
\\\hline
\end{tabular}
\end{center}
\end{table}

\begin{table}
\begin{center}
\caption{First 5 Newforms of Degrees 3, 4 and Weight 2\label{newforms3}}
\vspace{0.5ex}
\begin{tabular}{|l|l|l|l|}\hline
$k$ & $d$ & $N$ & reducible $p<10000$\\
\hline
2 & 3 & 41 & 17, 41 \\
2 & 3 & 53 & 13, 53 \\
2 & 3 & 61 & 61, 2087 \\
2 & 3 & 71 & 23, 31, 71, 479,  \\
&&&647, 1013, 3181\\
2 & 3 & 71 & 13, 71, 509, 3613 \\
\hline\end{tabular}
\hspace{2em}
\begin{tabular}{|l|l|l|l|}\hline
$k$ & $d$ & $N$ & reducible $p<10000$\\
\hline
2 & 4 & 47 & 47 \\
2 & 4 & 95 & 5, 19 \\
2 & 4 & 97 & 97 \\
2 & 4 & 109 & 109, 4513 \\
2 & 4 & 111 & 3, 37 \\
&&&\\
\hline\end{tabular}
\end{center}
\end{table}

\begin{table}
\begin{center}
\caption{First 5 Newforms of Degrees 2, 3 and Weight 4\label{newforms2w4}}
\vspace{0.5ex}
\begin{tabular}{|l|l|l|l|}\hline
$k$ & $d$ & $N$ & reducible $p<1000$\\
\hline
4 & 2 & 11 & 11 \\
4 & 2 & 13 & 13 \\
4 & 2 & 21 & 3, 7 \\
4 & 2 & 27 & 
{\tiny 3, 7, 13, 19, 31, 37, 43, 61, 67, 73, 79, 97, 103,}\\
&&&{\tiny 109, 127, 139, 151, 157, 163, 181, 193, 199, 211,}\\ 
&&&{\tiny 223, 229,241, 271, 277, 283, 307, 313, 331, 337,}\\
&&&{\tiny 349, 367, 373, 379, 397, 409, 421, 433, 439, 457,}\\
&&&{\tiny 463, 487, 499, 523, 541, 547, 571, 577, 601, 607,}\\
&&&{\tiny 613, 619, 631, 643, 661, 673, 691, 709, 727, 733,}\\
&&&{\tiny 739, 751, 757, 769, 787, 811, 823, 829, 853, 859,}\\ 
&&&{\tiny 877, 883, 907, 919, 937, 967, 991, 997} \\
&&&(has inner twists)\\
4 & 2 & 29 & 29 \\
\hline
\end{tabular}
\begin{tabular}{|l|l|l|l|}\hline
$k$ & $d$ & $N$ & reducible $p<1000$ \\\hline
4 & 3 & 17 & 17 \\
4 & 3 & 19 & 19 \\
4 & 3 & 35 & 5, 7 \\
4 & 3 & 39 & 3, 13 \\
4 & 3 & 41 & 41 \\
&&&\\&&&\\&&&\\&&&\\&&&\\&&&\\
&&&\\
&&&\\
\hline
\end{tabular}
\end{center}
\end{table}

\begin{table}
\begin{center}
\caption{Newforms on $\Gamma_0(389)$ of Weight $2$\label{newforms389}}
\vspace{0.5ex}
\begin{tabular}{|l|l|l|l|}\hline
$k$ & $d$ & $N$ & reducible $p<10000$ \\\hline
2 & 1 & 389 & none (degree 1 polynomials are all irreducible) \\
2 & 2 & 389 & {\tiny 5, 11, 59, 97, 157, 173, 223, 389, 653, 739, 859, 947,
1033, 1283, 1549, 1667, 2207, 2417, 2909, 3121, 4337,}\\
&&&{\tiny 5431, 5647, 5689, 5879, 6151, 6323, 6373, 6607, 6763, 7583, 7589, 8363, 9013, 9371, 9767} \\
2 & 3 & 389 & 7, 13, 389, 503, 1303, 1429, 1877, 5443 \\
2 & 6 & 389 & 19, 389\\
2 & 20 & 389 & 389 \\
\hline
\end{tabular}
\end{center}
\end{table}


\newpage
\begin{thebibliography}{XXXX}

\bibitem[BM03]{BM2003}
Baba, Srinath and Murty, Ram, \emph{Irreducibility of Hecke
Polynomials}, Math. Research Letters, 10(2003), no.5-6, pp.709-715.

\bibitem[FJ02]{fj}
D.~W. Farmer and K.~James, \emph{The irreducibility of some level 1 {H}ecke
  polynomials}, Math. Comp. \textbf{71} (2002), no.~239, 1263--1270.

\bibitem[FH91]{FH}
Fulton, William and Harris, Joe, \emph{Representation Theory, A
First Course}, Springer, 1991.

\bibitem[JO98]{JOno}
James, Kevin and Ono, Ken, \emph{A note on the Irreducibility of
Hecke Polynomials}, Journal of Number Theory 73, 1998, pp.527-532.

\bibitem[R80]{Rib80}
Ribet, Kenneth A., \emph{Twists of modular forms and endomorphisms
of abelian varieties.}, Math. Ann. 253 (1980), no. 1, 43--62.

\bibitem[R85]{Rib85}
Ribet, Kenneth A., \emph{On l-adic representations attached to
modular forms. II.}, Glasgow Math. J. 27 (1985), 185--194.

\bibitem[S81]{Serre}
Serre, Jean-Pierre,
\emph{Quelques applications du th\'eor\`eme de densit\'e de Chebotarev.}
Inst. Hautes Études Sci. Publ. Math. No. 54 (1981), 323--401.

\bibitem[SAGE]{sage}
Stein, William, \emph{Sage {M}athematics {S}oftware ({V}ersion 2.8.12)}, The
  SAGE~Group, 2007, {\tt http://www.sagemath.org}.

\end{thebibliography}
\end{document}